\documentclass{article}

\usepackage{graphicx}
\usepackage{caption}
\usepackage{subcaption}
\usepackage{float}
\usepackage{amsmath}
\usepackage{amsfonts}
\usepackage{amsthm}

\newtheorem{theorem}{Theorem}

\begin{document}

\begin{titlepage}
   \begin{center}

       \textbf{A generalization of the root function}

       \vspace{0.5cm}
       \textit{In memory of professor L\'aszl\'o Varga}
            
       \vspace{1.5cm}

	   Tam\'as D\'ozsa, \\
	   Ferenc Schipp

       \vspace{0.8cm}
            
       Department of Numerical Analysis\\
       E\"otv\"os Lor\'and University\\
       Budapest, Hungary\\
      
   \end{center}
   \vspace{0.8cm}
   
\textbf{Abstract.} We consider the interpretation and the numerical construction of the inverse branches of $n$ factor Blaschke-products on the disk $\mathbb{D}$ and show that these provide a generalization of the $n$-th root function. The inverse branches can be defined on pairwise disjoint regions, whose union provides the disk. An explicit formula can be given for the $n$ factor Blaschke-products on the torus, which can be used to provide the inverse branches on $\mathbb{T}$. The inverse branches can be thought of as the solutions $z=z_t(r) (0\le r\le 1)$ to the equation $B(z )=re^{it}$, where $B$ denotes an $n$ factor Blaschke-product. We show that starting from a known value $z_t(1)$, any $z_t(r)$ point of the solution trajectory can be reached in finite steps. The appropriate grouping of the trajectories leads to two natural interpretations of the inverse branches (see Figure \ref{fig:class_branch}). We introduce an algorithm which can be used to find the points of the trajectories.   

\vspace{6cm}

\small{EFOP-3.6.3-VEKOP-16-2017-00001: Talent Management in Autonomous Vehicle Control Technologies – The Project is supported by the Hungarian Government and co-financed by the European Social Fund.}

\end{titlepage}

%
%\keywords{Blaschke-functions, Blaschke-products, root-function, Newton-iteration}
%
%\mathclass{Primary 30J10, Secondary 26C15, 30J99, 65T99.}
%

%

%%%%%%%
% We shall fill out "???"
%%%%%%%%

%\commby{???}
%\recacc{???}{???}

\section{Introduction}
In this work we consider a generalization of the $n$-th root function $w\to \sqrt[n]{w}\ (w\in \Bbb C)$. Using polar coordinates 
$w=re^{it}\ (0\le r<\infty, t\in \Bbb I:=[-\pi,\pi))$,

\noindent we can express the branches of the root function as

\begin{equation}
\label{eq:n_root_branches}
\begin{split}
&\varphi_k(r,t):=r^{1/n} e^{i(t/n+2k\pi/n)}\\
&(r>0,  t\in \Bbb I, k=0,1,\dots,n-1).
\end{split}
\end{equation}

\begin{figure}[h!]
    \centering
    \begin{subfigure}{.5\textwidth}
        \centering
        \includegraphics[width=0.8\linewidth, trim={10cm 0 0 0},clip]{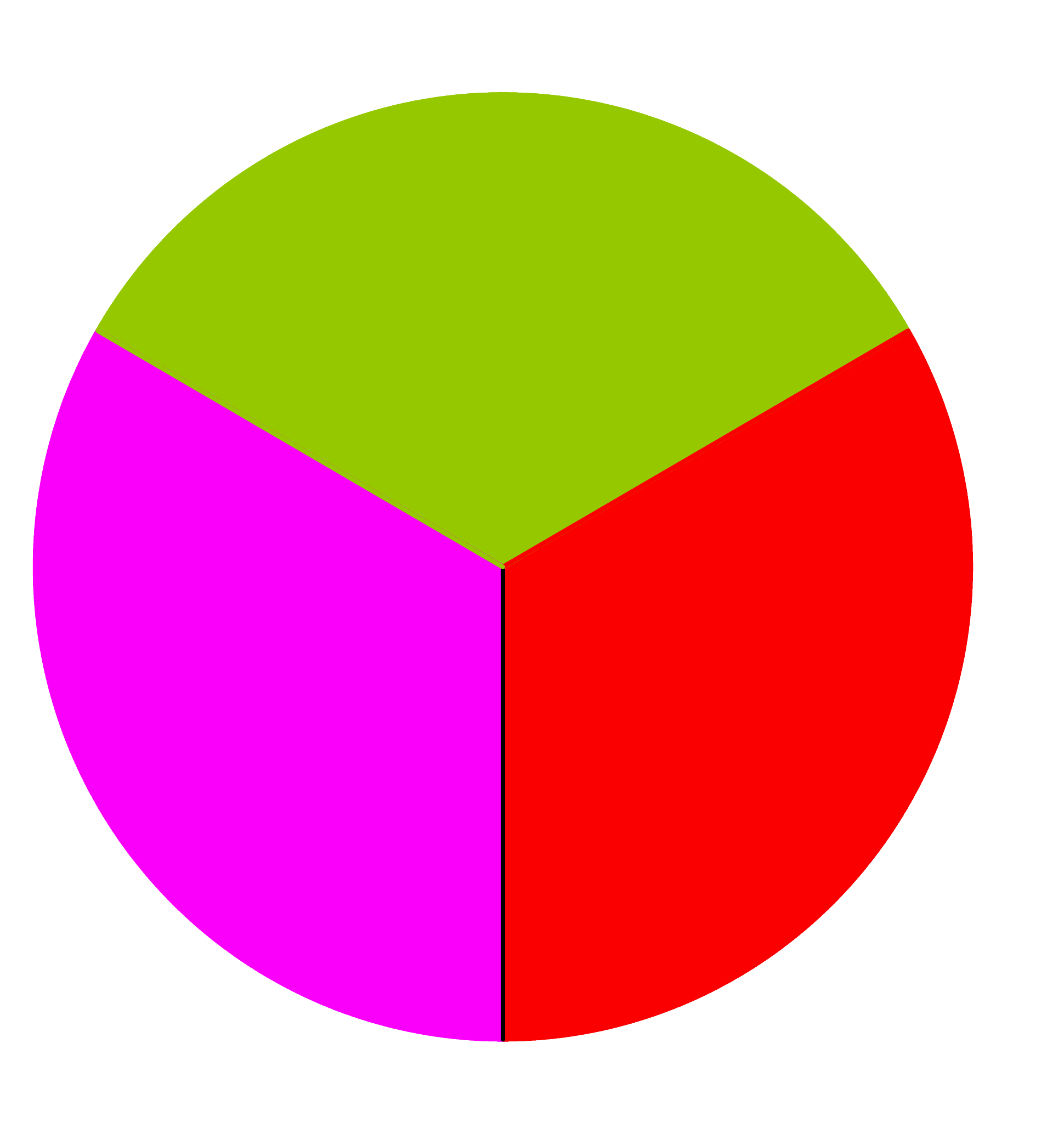}
        \caption{Ranges of $\varphi_k$}
        \label{fig:n_root_a}
    \end{subfigure}%
    \begin{subfigure}{.5\textwidth}
        \centering
        \includegraphics[width=0.9\linewidth, trim={0 0 0 10cm},clip]{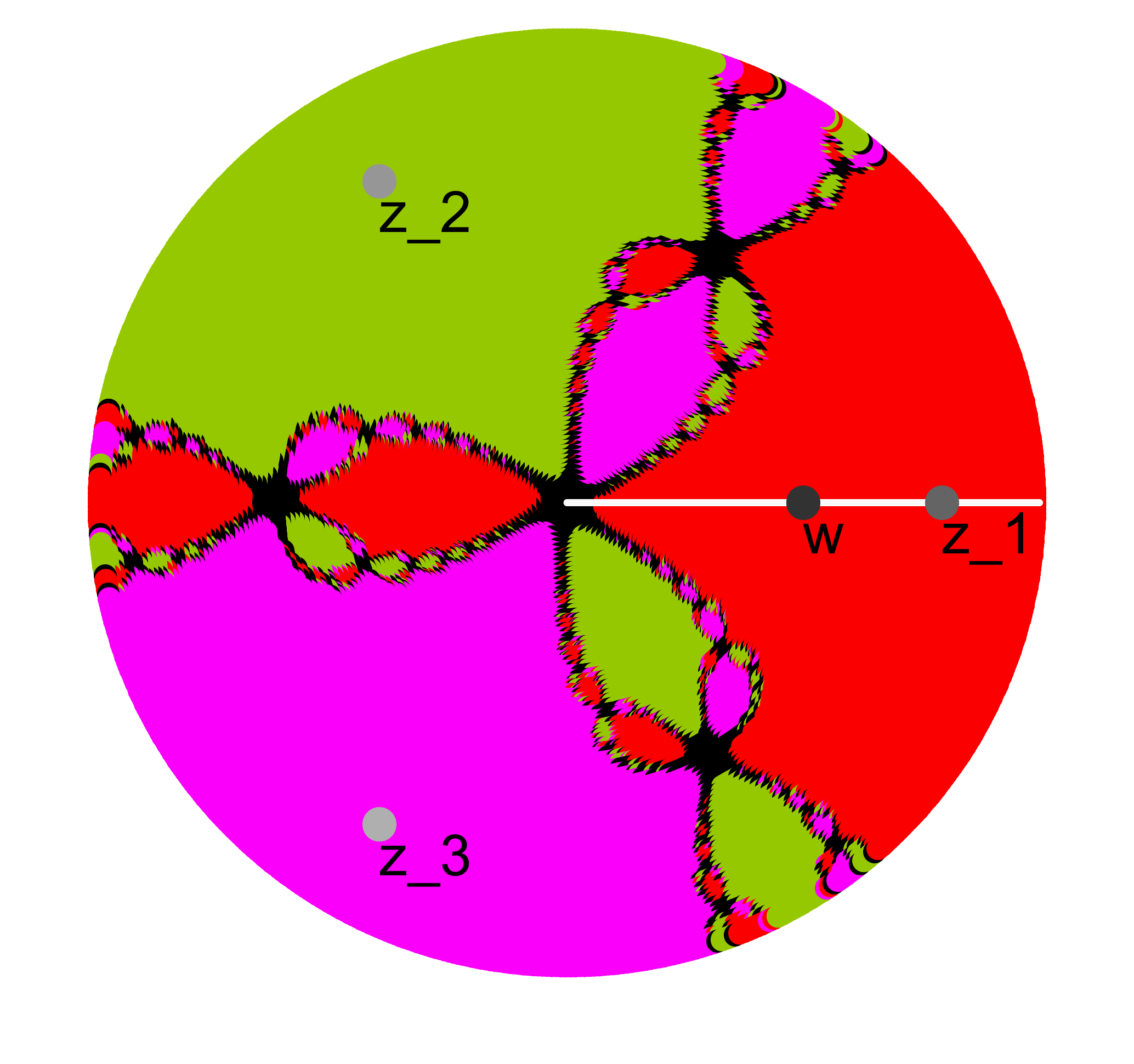}
        \caption{Initial points of Newton's method}
        \label{fig:n_root_b}
    \end{subfigure}
    \caption{Branches of the $n$-th root function}
\end{figure}

\noindent One can acquire all of the branches from a single one, using rotations of angle $2\pi/n$. The roots' absolute values $r^{1/n}$ and angles can be found with Newton's method and division respectively. Figure (\ref{fig:n_root_a}) shows the range of the branches of the $n$-th root function for $n=3$.

The solution of the equation $z^n=w \ (w\in\Bbb R)$ is equivalent to the $n$-th root function's value at $w$. In figure (\ref{fig:n_root_b}) we illustrate the roots of $w$ for $n=3$ on the unit disk. Each color represents a single root denoted by $z_1,z_2, z_3$. The roots can be identified as the limit of the Newton iteration

\begin{equation}
	\label{eq:Newton_it_4_root}
	v_{k+1}=v_k-\frac{g(v_k)}{g'(v_k)},\ g(v)=v^n-w\ \ (v_k\in\Bbb D, k=0,1,\dots).
\end{equation}

\noindent Any $v_0\in\Bbb D$ initial points for which the method diverges are colored black, the rest of the colors show the basins of attractions for each root $z_k$. The figure shows that choosing an initial point close to a particular root will result in Newton's method (quadratically) converging to that root.

In this paper we propose a generalization of the above problem by considering $n$ factor Blaschke-products instead of the $n$-th root function \cite{Bl}. Blaschke-functions have proven to be fundamental tools in the factorization of analytic functions and the construction of rational bases, while playing an important role for many applications in signal processing and control theory \cite{BSS, BSSA, BOSZA, HHW, HT, MoH}. The Blaschke-functions
\begin{equation}
	\label{eq:Blaschke_function}
   B_a(z):=\frac{z-a}{1-\overline a z}\ \ \left(z\in\Bbb C, a\in\Bbb D:=\{z\in\Bbb C:|z|<1\} \right)
 \end{equation}
 
 \noindent are bijections on the disk $\Bbb D$ and torus $\Bbb T:=\{z\in\Bbb C:|z|=1\}$. They can also be used to describe congruent transformations in the Poincar\'e model of the Bolyai--Lobachevskian geometry \cite{BOSZA, SCH}. The $n$-factor Blaschke-products

\begin{equation}
\label{eq:Blaschke_product}
\begin{split}
&\mathrm B(z):=\epsilon\prod_{k=1}^{n}\frac{z-a_k}{1-\overline a_k z }\\
(z\in\Bbb C,\epsilon\in &\Bbb T, A:=\{a_k:
1\le k\le n\}\subset\Bbb D)
\end{split}
\end{equation}

\noindent map $\Bbb D$, $\Bbb T$ and their complements into themselves. We refer to the set $A':=\{b\in\Bbb C:B'(b)=0\}$ as the set of {\it critical points} \cite{SS} and to the values $B(b), \ b \in A'$ as \textit{critical values}. If a line segment $S_t:=\{re^{it}:0\le r\le 1\}$ takes a critical value for some $r$, then $t$ is referred to as a critical parameter. The set of critical parameters will be denoted by $T_0$. We note that if $a_1=\dots=a_{n}=0,\epsilon=1$, then $\mathrm B(z)=z^n\ (z\in\Bbb C)$ therefore the inverse of $\mathrm B$ can be considered as a generalization of the $n$-th root function.

The Riesz-factorization provides the basis for several algorithms which deal with the rational represenation of signals \cite{CP, CS, TZW}. It can be considered as a generalization of polar coordinates for functions in $H^p(\Bbb D)$\ $(0<p\le\infty)$ Hardy-spaces. In this work instead of $H^p(\Bbb D)$, we consider the class of functions $\mathcal A_R$, which are analytic on $\Bbb D_R:=\{z\in\Bbb C:|z|<R\}$. For this class, if the roots of a function $f$ are known, its Riesz-factorization can be easily acquired. Namely, let $a_1,\dots, a_{n}$ denote the zeros of $f\in \mathcal A_R$. Then the factorization is given as

\begin{equation}
\label{eq:Riesz_fac}
\begin{split}
f(z)=&f_1(z)\prod_{k=1}^{n}(z-a_k)=\\
&=\mathrm B(z) f_1(z)\prod_{k=1}^{n}(1-\overline a_k z)=\mathrm B(z) S(z) \ (z\in \Bbb D_R),
\end{split}
\end{equation}
where the function $S\in \mathcal A_R$ has no zeros on $\Bbb D$. The existence of such a factorization, when $f\in H^p(\Bbb D)$ is however not easy to prove \cite{MoH}. Since $S(z)\ne 0\ (z\in\Bbb D)$ and $|\mathrm B(z)|=1\ (z\in\Bbb T)$, this factorization is often considered a generalization of polar coordinates.

It is well known that $\mathrm B:\Bbb T\to\Bbb T$ provides an $n$-fold mapping of $\Bbb T$, moreover there exists a $\beta=\beta_n:\Bbb R\to\Bbb R$ strictly increasing function $\beta(t+2\pi)=\beta(t)+2\pi\ (t\in\Bbb R)$ for which

\begin{equation}
\label{eq:perem_Blaschkek}
\mathrm B(e^{it})=e^{in\beta(t)}\ \
(t\in\Bbb R).
\end{equation}

\noindent
The function $\beta$ can be expressed in an explicit form if the parameters of the Blaschke-product are known. As a result, we can easily produce the inverse of a Blaschke-product on the torus \cite{SCH}. Namely, for the function

$$
\tau_k(t):=\beta^{-1}(t/n+2k\pi/n)\ \ (t\in\Bbb R, k=0,\dots, n-1),
$$
the following holds
$$
\mathrm B(e^{i\tau_k(t)})=e^{in\beta(\beta^{-1}(t/n+2k\pi/n))}=e^{it}.
$$

\noindent
From the above, we get that $\varphi_k(1,t):=e^{i\tau_k(t)}\ (t\in\Bbb R)$ are the inverse branches of $\mathrm B$ on the torus.

In general, finding the inverse branches can be reduced to finding the roots of an $n$-th degree polynomial. Unfortunately, when considering this form of the problem, accurate geometric description of the results seems impossible. In the case of two factor Blaschke-products, the inverse branches can be expressed explicitly using the complex root function \cite{DS}. These results provide a basis for the description of the inverse branches for $n$ factor Blaschke-products.

In this paper, we consider a description of the inverse branches and propose a numerical algorithm to produce them. We plan to extend these results to functions in $\mathcal A_R$ in a future work, hence providing a numerical method to produce their Riesz-factorization. Notice, that when considering the inverse branches of $\mathrm B$, for each $w=re^{i t}\in\overline{\Bbb D}$, the equation $\mathrm B(z)=w$ is equavilent to finding the roots of an $n$ degree polynomial, therefore (considering multiplicities as well) has exactly $n$ solutions. When grouping these solutions into $n$ classes, there are several ways we can define (continous) inverse branches. In section 2, we shall prove the following theorem about the branches:

\begin{theorem}
For all
 $ t\in \Bbb T, t\notin T_0$ the equation $\mathrm B(z)=re^{it}$ has a single, unique (continuous in $r$) solution $\varphi_k(r,t)$, for which
 \begin{equation}
 \begin{split}
 &\mathrm B(\varphi_k(r,t))=re^{it}, \varphi_k(1,t)=e^{i\tau_k(t)}\\
 &(0\le r\le 1,  k=0,\dots,n-1),
 \end{split}
 \end{equation}
 furthermore these solution curves are disjoint for $r>0$ and their union provides the closed disk.
\end{theorem}

\begin{figure}[H]
    \centering
    \begin{subfigure}{.5\textwidth}
        \centering
        \includegraphics[width=0.9\linewidth, trim={0 0 0 10cm},clip]{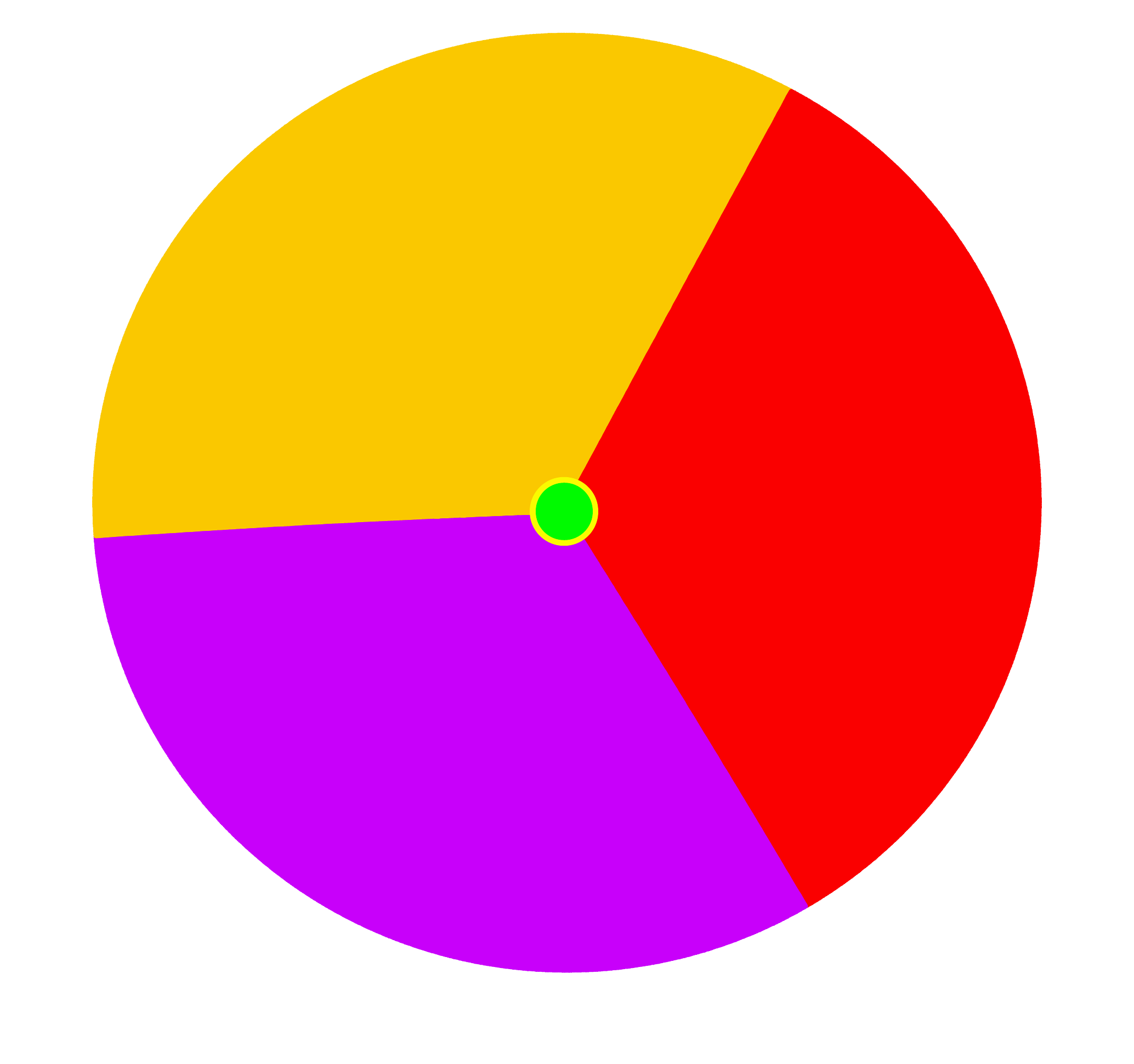}
        \caption{Inverse branches of \\ $a_1=\dots=a_3=0$}
        \label{fig:n_root_branches}
    \end{subfigure}%
    \begin{subfigure}{.5\textwidth}
        \centering
        \includegraphics[width=0.9\linewidth, trim={0 0 0 10cm},clip]{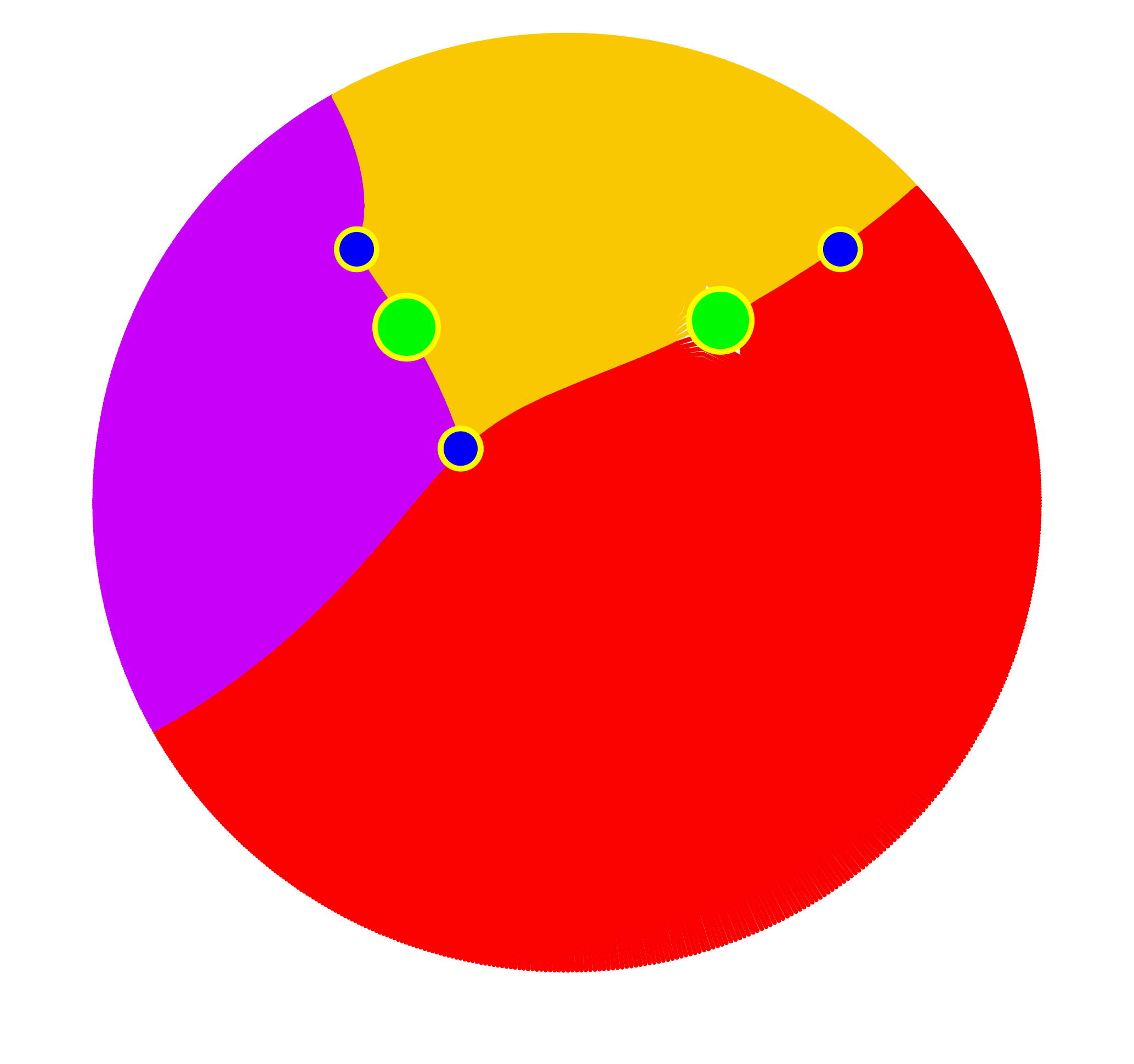}
        \caption{Theorem 1. based classification of branches}
        \label{fig:branches_theorem_one}
    \end{subfigure}
        \begin{subfigure}{.5\textwidth}
        \centering
        \includegraphics[width=0.9\linewidth, trim={0 0 0 10cm},clip]{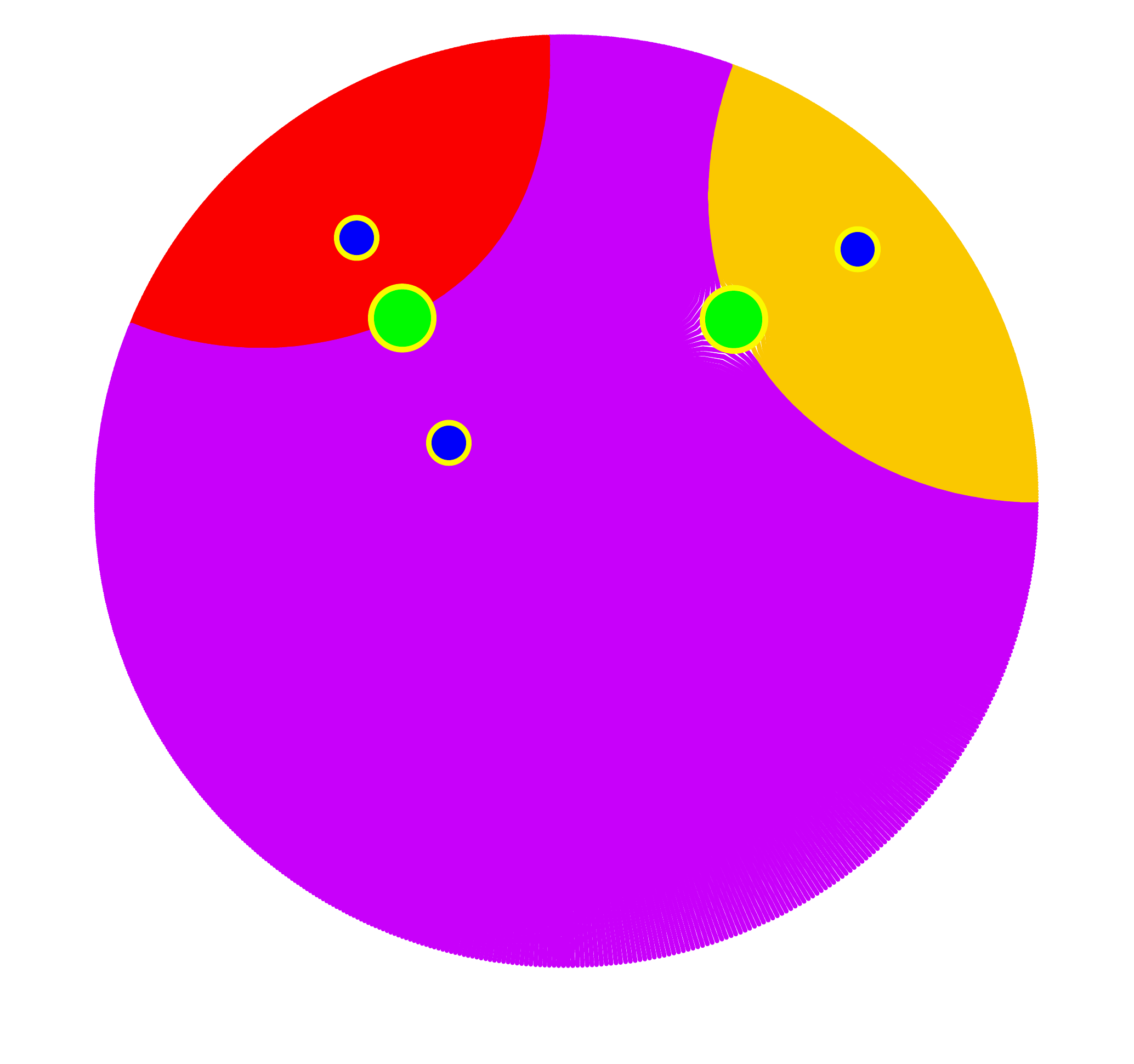}
        \caption{Common zero based classification of branches}
        \label{fig:branches_common_zero}
    \end{subfigure}
    \caption{Classification of inverse branches}
    \label{fig:class_branch}
\end{figure}

There are other ways of defining the inverse lines of the generalized root function as well. For example one could consider all trajectories ending in the same zero to be an inverse branch. Figure (\ref{fig:class_branch}) depicts the inverse branches for a $3$ pole Blaschke-product. In figures (\ref{fig:branches_theorem_one}) and (\ref{fig:branches_common_zero}), small circles show the zeros of $\mathrm B$, while the larger circles refer to the critical points. Figure (\ref{fig:n_root_branches}) shows the branches of the function $w\to\sqrt[3]{w}$, in which case $a_1=\dots=a_3=0$, or in other words $0$ is a root with a multiplicity of $3$. Figure (\ref{fig:branches_theorem_one}) illustrates the inverse branches as defined in Theorem 1, while in (\ref{fig:branches_common_zero}) a branch is made up of the trajectories ending in the same zero points. The dots show the positions of the inverse poles and critical points.

\section{Differential equation form of the problem}

The Blaschke-product $\mathrm  B$ has $2n-2$ critical points which fall in the hyperbolic convex hull of the zeros \cite{SS}. All other critical points fall outside $\Bbb D$ and are the reflections of the ones in $\Bbb D$ with respect to $\Bbb T$ \cite{SS}. Let $R:=\min\{1/|a|: a\in A\}$ and set of critical values by
$$
K:=\mathrm B(A'):=\{\mathrm B(z) :z\in A'\}.
$$

\noindent

We show that if $w=re^{it_0}\notin K\ (r\in [r_1,r_2])$ and $\mathrm B(z_0)=w_0$, then the implicit equation
\begin{equation*}
\mathrm B(\varphi(r))=w=re^{it_0}, \varphi(r_0)=z_0\ \ \ (r\in [r_1,r_2]\subset I:=(0,R))
\end{equation*}
is equivalent to an initial value problem corresponding to a first order differential equation \cite{A}. Since for all $r\in [r_1,r_2]$, $w=re^{it_0}\notin K$, the function $\mathrm B(z)=w$ has a locally unique inverse, piecing these together yields the differentiable solution $\varphi(r)=\varphi(r,t_0)\ (r\in [r_1,r_2])$ of the implicit equation, for which
\begin{equation}
\label{eq:implicit_alak}
\mathrm B(\varphi(r))=re^{it_0}, \  \varphi'(r)\mathrm B'(\varphi(r))=e^{it_0}\ \ (r\in [r_1,r_2]).
\end{equation}
We note that if the roots of $\mathrm B$ have a multiplicity of $1$, then the choice $r_1=0$ is also possible.

Let us introduce the function
$$
f(z):=\frac {\mathrm B(z)}{\mathrm B'(z)}\ \ (z\in \Omega:=\Bbb D_R\setminus (A\cup A')).
$$
We show that if $z_0\in\Omega$, then the problem (\ref{eq:implicit_alak}) is equivalent to the initial value problem for the function $\chi: I \to \Omega, \ I \subset (0, \infty)$:
\begin{equation}
\label{eq:diffegyenlet_alak}
\chi'(r)=\frac{f(\chi(r))}r, \ \ \chi(r_0)=z_0\ \
(r_0\in I,z_0\in \Omega)
\end{equation}
where $ \mathrm B(z_0)=w_0=r_0e^{it_0}\notin K$. Indeed, by (\ref{eq:implicit_alak})
$$
\varphi'(r) =
\frac{e^{it_0}}{\mathrm B'(\varphi(r))}=\frac{\mathrm B(\varphi(r))}{r\mathrm B'(\varphi(r))},
$$
therefore $\varphi$ is a solution of (\ref{eq:diffegyenlet_alak}). Conversely, for any $\chi$ solution of (\ref{eq:diffegyenlet_alak})
$$
  \frac d{dr}\log r=\frac 1 r=\chi'(r)\frac{\mathrm B'(\chi(r))}{\mathrm B(\chi(r))}
 =\frac d{dr}\log \mathrm B(\chi(r))
$$
holds, therefore $\mathrm B(\chi(r))=cr \ (r\in (r_1,r_2))$, where by the initial values, for the constant $c$ we get
$$
r_0 e^{it_0}=\mathrm B(z_0)=\mathrm B(\chi(r_0))=c \cdot r_0,
$$
from which $c=e^{it_0}$ follows.

The right side of (\ref{eq:diffegyenlet_alak}) fulfills the usual conditions for the existance and uniqueness of the solution, therefore for any $(r_0,z_0)\in I\times\Omega$, the differential equation (\ref{eq:diffegyenlet_alak}) has a unique extended solution $\chi_{z_0}(r)\ (r\in J\subset I)$ to the boundary \cite{A} (pp. 50-57). The ranges of these solutions will be referred to as trajectories henceforth. From the existence and uniqueness theorem, it follows that two trajectories are either disjoint or are the same. Furthermore, for any point in $\Omega$ there is a trajectory which passes through it, therefore taking the union of these trajectories yields $\Omega$.

In order to describe the maximal solutions, let us introduce the partitioning of the line segments $S_t$ induced by the finite set $K$. For 
$$
S_t\cap K:=\{ \rho_je^{it}:j=1,\dots, m_t\}
$$
let $\rho_0:=0, \rho_{m_t+1}:=R$, furthermore let
\begin{equation*}
\begin{split}
&\mathcal J_t:=\{J_j:=(\rho_{j-1}e^{it},\rho_{j}e^{it}):j=1,\dots, m_t+1\},\\
&\mathcal J_t:=\{(0,R)\},\ \ \text{if}\ \ S_t\cap K=\emptyset.
\end{split}
\end{equation*}
Denote line segments $J=(\rho e^{it},\sigma e^{it})\in \mathcal J_t$ by $\widehat J:=(\rho ,\sigma )$. Then the solutions $\varphi=\varphi_J \ (J\in\mathcal J_t)$ of the implicit equation
$$
\mathrm B(\varphi(r))=e^{it}r\ \ (r\in \widehat J, J\in\mathcal J_t, t\in \Bbb I)
$$
match the maximal solutions of (\ref{eq:diffegyenlet_alak}).

Now we can describe the inverse branches in the following manner. For any $w=re^{it_0}, t_0\notin T_0 $ take the initial values $z_{0k}=e^{i\tau_k(t_0)}$ $(k=0,\dots,n-1)$. We note that for these values $\mathrm B(z_{0k})=e^{it_0}$ holds. Using these initial values, let us consider the maximal solutions of (\ref{eq:diffegyenlet_alak}) $\varphi_k=\chi_{z_{0k}}$. If the zeros have a multiplicity of one and $t_0\notin T_0$, the domain of these solutions is $[0,R)$, furthermore the solutions match the trajectories from Theorem 1.

\begin{figure}[h!]
    \centering
    \begin{subfigure}{.5\textwidth}
        \centering
        \includegraphics[width=0.9\linewidth, trim={0 0 0 10cm},clip]{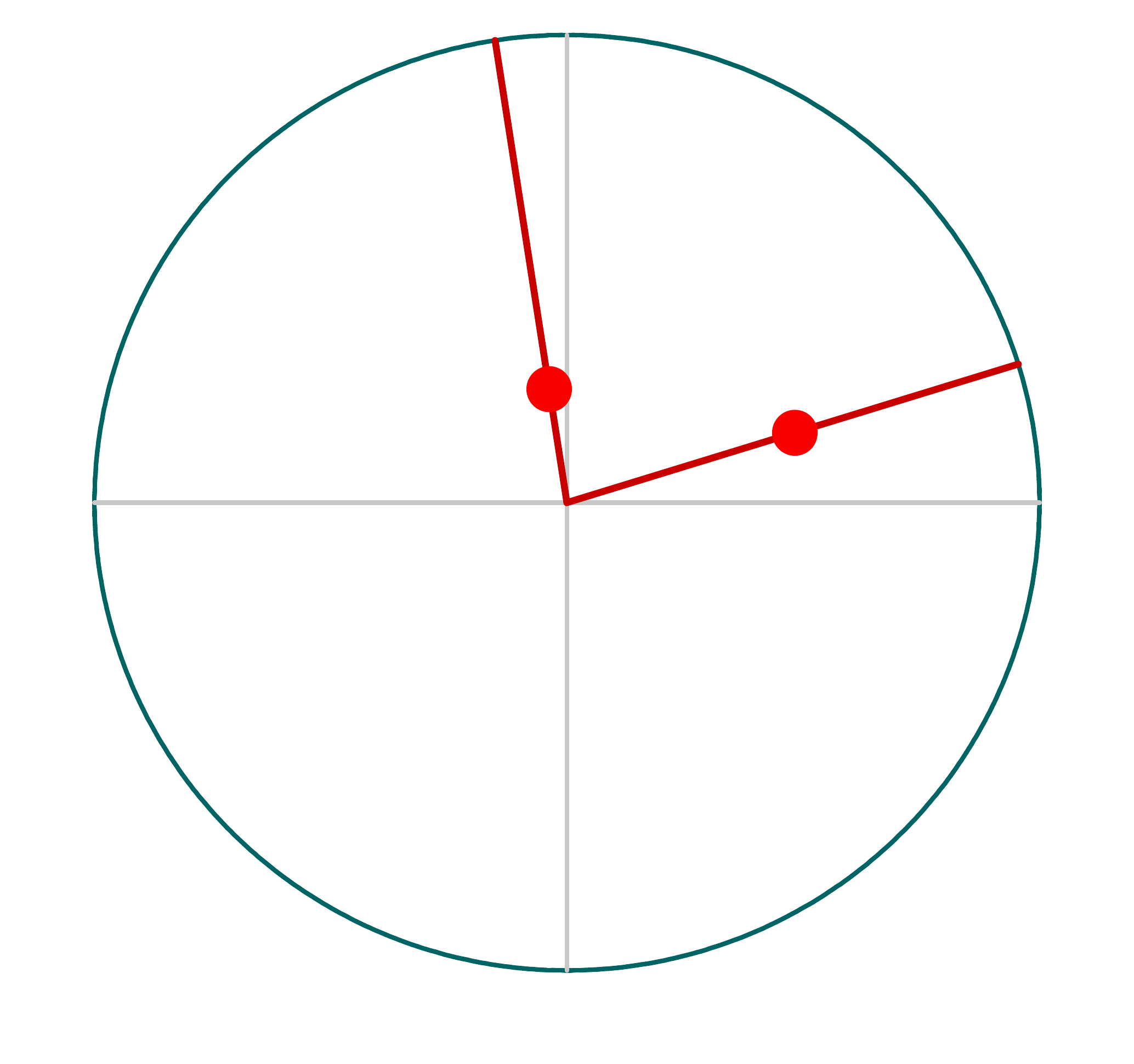}
        \caption{Critical lines $S_t$}
        \label{fig:crit_line_ims}
    \end{subfigure}%
    \begin{subfigure}{.5\textwidth}
        \centering
        \includegraphics[width=0.9\linewidth, trim={0 0 0 9cm},clip]{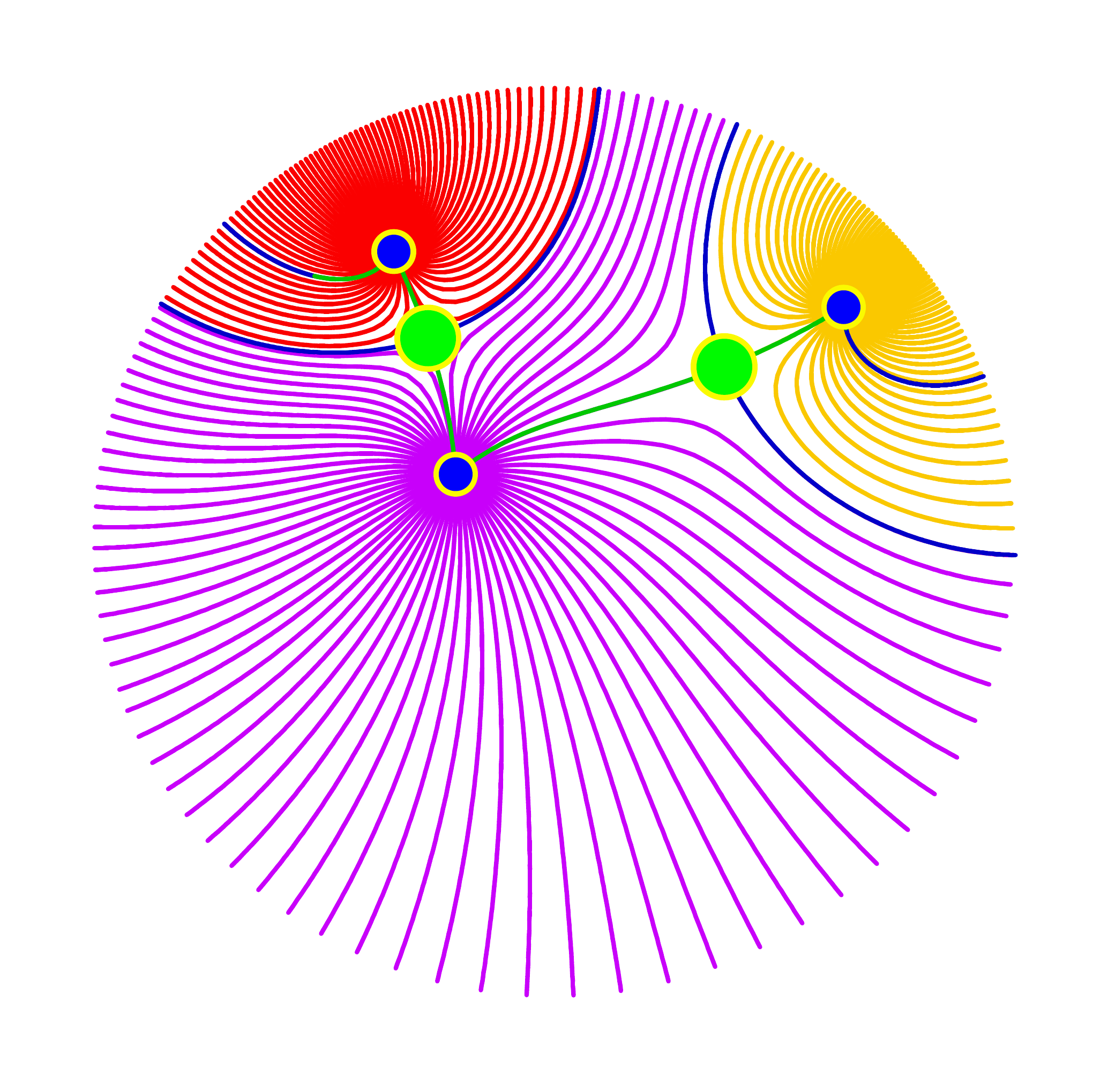}
        \caption{Inverse images of $S_t$ and the trajectories}
        \label{fig:crit_line_invs}
    \end{subfigure}
    \caption{Critical lines and their inverse images}
    \label{fig:crit_lines}
\end{figure}

Figure (\ref{fig:crit_lines}) illustrates the partitioning of the trajectories. Figure (\ref{fig:crit_line_ims}) shows the set $K=\mathrm B(A')$ (for $|K| = 2$) and the partitioning of the critical line segments $S_t\ (t\in T_0)$, which in this case consists of $2$ segments. The inverse images of these are illustrated in figure (\ref{fig:crit_line_invs}).

Taking nested circles instead of line segments, we get the following theorem.

\begin{theorem}
\it Suppose that for the circle $C_r:=\{re^{it}: t\in\Bbb I\}\ (0<r\le 1)$,
 $ C_r\cap K=\emptyset$ holds. Let $w_0=re^{it_0}\in C_r$ and denote by $z_{0j}$  the inverse images of $w_0$:
 $B(z_{0j})=w_0\ (j=0,\dots,n-1)$. Then, the smooth solutions of the implicit equation
 \begin{equation}
 \mathrm B(\psi_j(t))=re^{it}\ (t\in \Bbb I),\ \psi_j(t_0)=z_{0j}
 \end{equation}
exist uniquely.
\end{theorem}

\begin{figure}[h!]
    \centering
    \begin{subfigure}{.5\textwidth}
        \centering
        \includegraphics[width=0.9\linewidth, trim={0 0 0 10cm},clip]{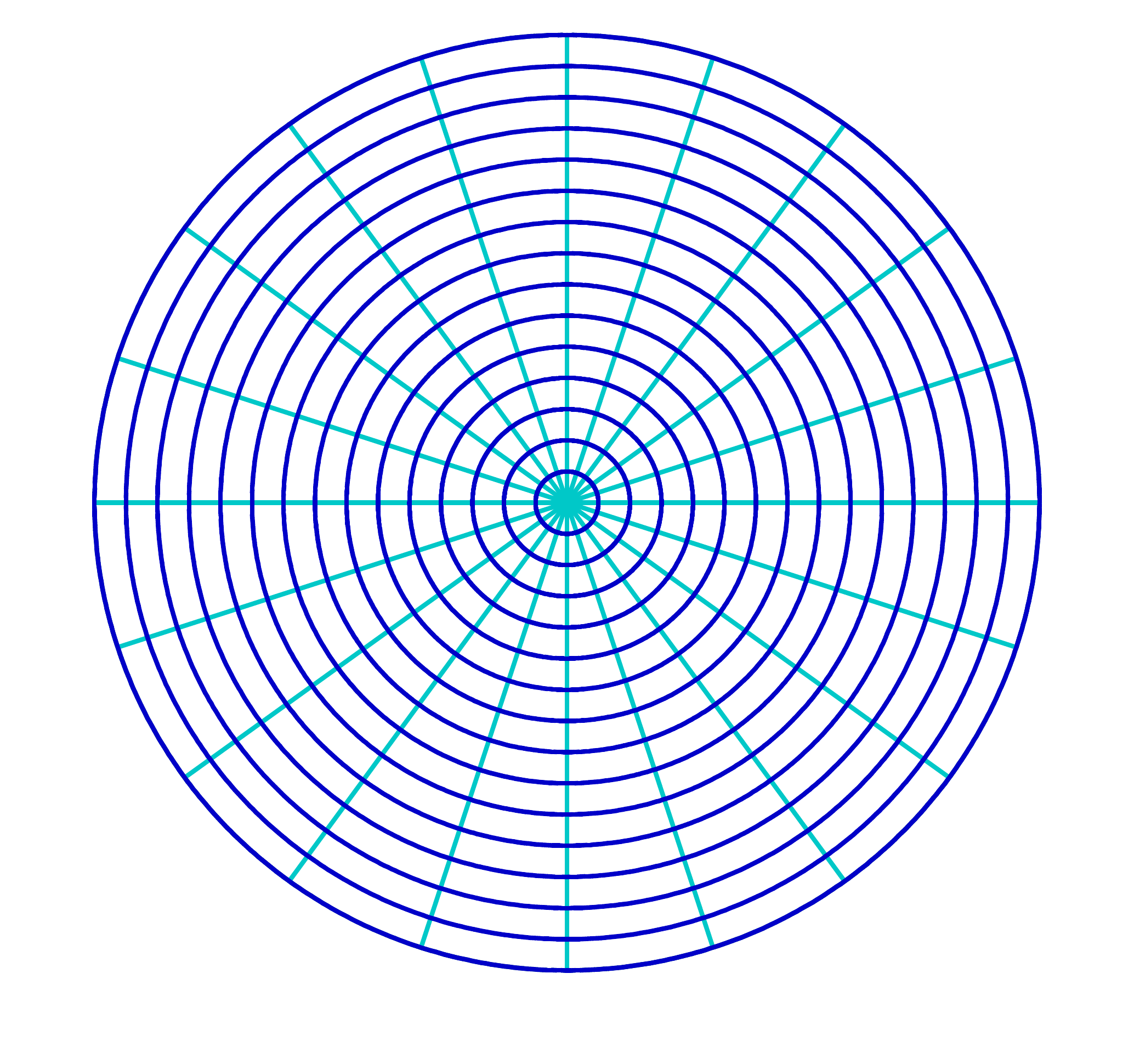}
        \caption{Nested circles}
        \label{fig:nested_circle_ims}
    \end{subfigure}%
    \begin{subfigure}{.5\textwidth}
        \centering
        \includegraphics[width=0.9\linewidth, trim={0 0 0 10cm},clip]{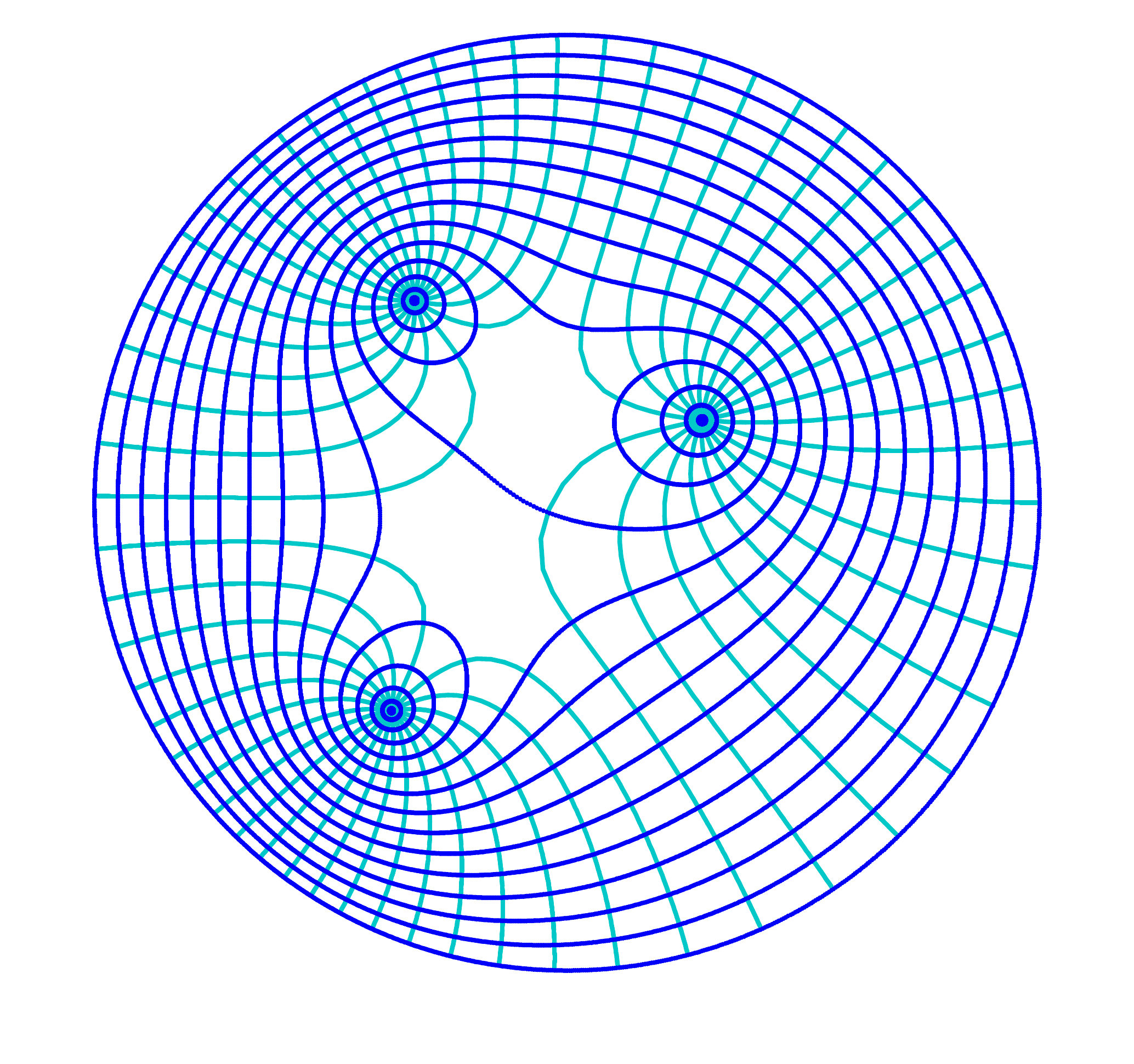}
        \caption{Inverse images of nested circles}
        \label{fig:nested_circle_invs}
    \end{subfigure}
\end{figure}

Figures (\ref{fig:nested_circle_ims}) and (\ref{fig:nested_circle_invs}) show the nested circles and the corresponding inverse images.

\section{Inverse algorithm}

In this section we consider a numerical solution to the problem outlined above. The problem of finding the inverse of a Blaschke-product can be posed in three equivalent forms. Namely, the inverse images are solutions to the following implicit equation, differential equation and line integral problems.  

\begin{equation*}
\begin{split}
&i)\ \ \ \mathrm B(z)=w=re^{it}, B(z_0)=w_0\ \ (z=\varphi(r))\\
&ii)\ \ \chi'(r)=f(\chi(r))/r, \chi(r_0)=z_0\ \ \ (f(z)=\mathrm B(z)/ \mathrm B'(z))\\
&iii)\ F(z) = \int_{z_0 }^z\frac{\mathrm B'(\zeta)}{\mathrm B(\zeta)}\, d\zeta
\end{split}
\end{equation*}
Any of these forms can be solved with various efficient numerical algorithms \cite{Hen}. In this section we introduce an algorithm based on Newton's method which solves the implicit equation form of the problem. We generalize the findings of the previous section and examine the solutions of 
\begin{equation}
\label{eq:impl_num}
  \mathrm B(z)=w\in\Gamma,
\end{equation}
where $\Gamma\subset\Bbb D_R$, $\Gamma\cap K=\emptyset$ is a simple smooth curve. We show that starting from the initial solution $\mathrm B(z_0)=w_0$, one can acquire a finite number of  $(z_k,w_k)$ $(k=1,\dots,N)$ solutions for which $(z,w)=(z_N,w_N)$ and $z_{k+1}$ is a limit point of a Newton-iteration starting from $z_{k}$. These results are a consequence of the following theorem.

\begin{theorem}
\it  There exist numbers $r_0>0, r_1>0$, such that for any $\mathrm B(z_0)=w_0, |z_0|\le 1$ initial solution and any $w\in K_{r_1}(w_0)\cap \Gamma$ a unique $z\in K_{r_0}(z_0),|z|\le 1, \mathrm B(z)=w$ solution exists of (\ref{eq:impl_num}), which is the limit of the quadratically convergent sequence 
$$
z_{j+1}=g(z_j),\ \ (j\in\Bbb N)
$$
where
\begin{equation}
\label{eq:Newton_it}
\begin{split}
 &g(z):=g_{w}(z)=(z-z_0)-\frac{\mathrm B(z)-w}{\mathrm B'(z)}\ \ (z\in K_{r_0}(z_0))\\
 &z=\lim_{j\to\infty} z_j.
 \end{split}
\end{equation}
\end{theorem}

\proof We show the existence of $r>0$ (independent of $z_0$), such that the mapping $g$ is contractive on the disk $Z:=\{z\in \Bbb C:|z-z_0|\le r\}$:
$$
g:Z\to Z,\ \  |g(z_1)-g(z_2)|\le \frac 12|z_1-z_2|\ \ (z_1,z_2\in Z).
$$

The function $g$ is twice differentiable on $\overline{\Bbb D}_R$ and its derivatives satisfy
\begin{equation*}
\begin{split}
& g'(z)=
\frac{(\mathrm B(z)-w)\mathrm B"(z)}{(\mathrm B'(z))^2},\\
& g"(z)=
\frac{\mathrm B"(z)}{\mathrm B'(z)}+\frac{(\mathrm B(z)-w)\mathrm B^{(3)}(z)}{(\mathrm B'(z))^2}-
2\frac{(\mathrm B(z)-w)(\mathrm B"(z))^2}{(\mathrm B'(z))^3}
\end{split}.
\end{equation*}

First we prove that there exists $\rho_0>0$, independent of $z_0$, so
\begin{equation}
\label{eq:mean_val}
w_0=\mathrm B(z_0)\in \Gamma \ \ \Rightarrow\ \ \ \forall \kappa\in K:\ |\kappa-z_0|\ge
\rho_0.
\end{equation}
Indeed, by the mean value theorem $\mathrm B'(\kappa)=0, M:=\sup_{|z|\le R}|\mathrm B"(z)|<\infty$
we get
$$
0<d:=\min_{\kappa\in K, w\in\Gamma} |\kappa-w|
\le|\mathrm B(\kappa)-\mathrm B(z_0)|\le
M|\kappa-z_0|^2,
$$
from which (\ref{eq:mean_val}) holds with the constant $\rho_0:=\sqrt{d/M}$. Let
$$
 H:=\overline{\Bbb D}_R
\setminus\bigcup_{\kappa\in K} K_{\rho_1}(\kappa)\ (\rho_1:=\rho_0/2).
$$
The functions $g',g"$ are bounded on $H$, furthermore $\overline K_{\rho_1}(z_0)\subset H$ holds. Let
$$
m:=\min_{z\in H}|\mathrm B'(z)|,\ \  M_1:=\max_{|z|\le R}|g"(z)|.
$$
Then by $M_1<\infty,m>0$ and the mean value theorem
\begin{equation}
\label{eq:g_der_felsobecsles}
\begin{split}
&|g'(z)|\le |g'(z)-g'(z_0)|+|g'(z_0)|\le\\
 &\le M_1|z-z_0|+\frac {M}{m^2}|w_0-w|<\frac 12\
 \left(z\in K_{\rho_1}(z_0) \right)
 \end{split}
  \end{equation}
  holds if, for example
  $$
  |w-w_0|\le \frac {m^2}{4M},\ |z-z_0|\le \frac 1{2M_1}.
  $$
  By the above, for $r=\min\{\rho_1, 1/(2M_1)\}$ the mapping $g $ is contractive:
  $$
  |g(z_1)-g(z_2)|\le \frac 12|z_1-z_2|\ \ \left(z_1,z_2\in \overline K_r(z_0) \right)\ .
  $$
  Since
  \begin{equation}
  \begin{split}
  &|g(z)-z_0|\le |g(z)-g(z_0)|+|g(z_0)-z_0|\le\\
  &\le \frac 12 |z-z_0|+\frac{|w_0-w|}m,
  \end{split}
  \end{equation}
  choosing the parameters $r,r_0,r_1$ according to
  \begin{equation}
  \begin{split}
  &r=r_0\le \min\{\rho_1,1/(2M_1)\},\ K_r(z_0)\subset K_R(0),\\
  &|w-w_0|\le r_1:=\min\{mr/2,m/(2M)\}
  \end{split}
  \end{equation}
  the mapping $g:Z:=\overline K_r(z_0)\to Z$ is indeed contractive, with a contraction coefficient of $1/2$.

  By the fixed-point theorem the limit  
  $$
  z:=\lim_{j\to \infty}z_j\in Z
  $$
  exists and
  $$
  g(z)=z\ \ \Leftrightarrow\ \ \mathrm B(z)=w.
  $$
  Finally, $|z|\le 1$ holds, since $w\in \Gamma\subset \overline{\Bbb D}$ and $\mathrm B$ maps $\overline{\Bbb D}$ onto itself. Since $g(z)=z$, the derivative $g'(z)=0$ and by     the mean value theorem
  $$
  |z_{j+1}-z|=|g(z_j)-g(z)|\le\max_{\zeta \in Z}|g"(\zeta)||z_j-z|^2/2\le M_1|z_j-z|^2/2,
  $$
  meaning that the method is quadratically convergent, which concludes our proof. \qed
                         
  We would like to draw attention to the fact, that in applications calculating the derivative function $B'(z)$ in (\ref{eq:Newton_it}) may pose some difficulty. In \cite{SAK} the authors describe a method with which the derivatives of a polynomial can be calculated from its roots and the function value. The same idea can be easily extended for Blaschke-products and a simple formula for $B'(z)$ expressed with the function value and the poles can be given:
\begin{equation}
\label{eq:der_blasch}
B'(z) = B(z)\cdot \left[ \sum_{i=1}^{n} \frac{1}{z-a_i} + \sum_{i=1}^{n} \frac{\overline{a_i}}{1-\overline{a_i}z} \right].
\end{equation}
Suppose we have $w_0, z_{0i}$ such that $B(z_{0i}) = w_0, \ (i=1,\ldots,n)$. Then, $z_{0i}$ are the roots of the polynomial
$$
P(z) = \prod_{i=1}^{n} (z-a_i) - w_0 \prod_{i=1}^{n}(1-\overline{a_i}z),
$$
furthermore $P(z_{0i}) = B(z_{0i}) - w_0 = 0$. In a future work we hope to develop a method to calculate the derivative values from only the known $z_{0i}, \ (i=1,\ldots,n)$ based on the above observations.

\end{document}